\documentstyle[12pt]{article}

\begin{document}
\vspace{1cm}
\qquad \qquad \qquad \qquad \qquad  \qquad \qquad \qquad
{\bf In Memory of Victor Nikolaevich Popov}
\vspace{1cm}

\begin{center}

{\large {\bf GAUSS DECOMPOSITION}}\\
{\large {\bf FOR QUANTUM GROUPS AND SUPERGROUPS}} \\[.5cm]

{\bf E.V.Damaskinsky} \footnote
 {  Supported by Russian Foundation for
 Fundamental Research, Grant N 94-01-01157-a.}  \\
 Defense Constructing Engineering Institute,
 Zacharievskaya st 22, \\ 191194, St Petersburg, Russia \\[.5cm]

{\bf P.P.Kulish} \footnote{ E-mail address:
kulish@pdmi.ras.ru \quad  kulish@evalvx.ific.uv.es}\\
St.Petersburg Branch of Mathematical Institute of
Russian Academy of Science \\
Fontanka 27, 191011, St.Petersburg, Russia \\
Department of Theoretical Physics, Valencia University, Spain \\
Erwin Schr\"odinger International Institute for Mathematical
Physics (ESI)\\
Pasteurgasse 6/7, A-1090 Wien, Austria \\[.5cm]

{\bf M.A.Sokolov} \footnote{ E-mail address:
 sokolov@pmash.spb.su} \\
 St Petersburg Institute of Mashine Building, Poliustrovskii pr 14, \\
 195108, St Petersburg, Russia \\[.5cm]

  {\bf Abstract.}
\end{center}
The Gauss decompositions of the quantum groups, related to classical
Lie groups and supergroups are considered by the elementary algebraic and
$R$-matrix methods. The commutation relations between new basis generators
(which are introduced by the decomposition) are described in some details.
It is shown that the reduction of the independent generator number
in the new basis
to the dimension of related classical (super) group is possible. The classical
expression for (super) determinant through the Gauss decomposition
generators is not changed in the deformed case. The symplectic quantum
group $Sp_q(2)$ and supergroups $GL_q(1|1), GL_q(2|1)$ are considered as
examples.

\newpage
{\bf 1. Introduction.} The development of the quantum inverse
scattering method (QISM) [1] intended for investigation of the
integrable models of the quantum field theory and
statistical physics naturally gives rise to some
interesting algebraic constructions. Their investigation allows to select
a special class of Hopf algebras now known as quantum groups and
quantum algebras [2,3]. In many respects these algebraic structures
resemble very much those of the standard Lie groups and Lie algebras, that is
probably one of the reasons
explaining the appearance of quite a few number of articles
dealing with these objects. The nice $R$-matrix formulation of the quantum
group theory [4], based on the fundamental relation of
the QISM (the FRT-relation) has given an
additional impulse for these investigations. The FRT-relation
\begin{equation}
\label{1}
RT_1T_2=T_2T_1R,
\end{equation}
may be considered as a condensed matrix form of homogeneous quadratic
(commutation) relations between $n^2$ generators $t_{ij},\,\,
(i,j=1,2,...,n)$ of the associative algebra $F_q$.  The entries of the
number $n^2 \times n^2$ matrix $R$ play the role of the Lie algebra
structure constants. We are using the QISM standard
notations $T_1=T\otimes I,\em
T_2=I\otimes T$, where  $T=(t_{ij })$ is the $n\times n$-matrix of generators
for the associative algebra $F_q\,\,$ ($q$-{\bf matrix}). Usually it is
supposed
that the $R$-matrix is a solution of the Yang-Baxter equation [4].
The $R$-matrices related to the classical simple Lie algebras from
the Cartan list are well known [3-5]. The matrix form of the
FRT-relation (\ref{1}) allows to introduce without difficulty
into the associative algebra $F_q$, together with multiplication, additional
maps (coproduct, counit and antipod (provided that the inverse matrix
$T^{-1}$ exists)) and convert $F_q$ into the Hopf algebra i.e.
 the {\bf quantum group} [4].  Rather straightforward (and sometimes
misleading) analogy of the $q$-matrix $T$ with the standard matrices
as elements of the classical matrix
groups nevertheless allows one to consider many questions
of the higher algebra without difficulties operating with quantum analogs
of the well known
classical objects such as determinants, minors, linear spaces,
symmetric and external products etc.[4, 6-11].

In this paper we consider the Gauss decomposition of the $q$-matrix
$T=T_LT_DT_U$ on strictly lower and upper triangular matrices (with
unities on the diagonal) and the diagonal matrix
$T_D:= A$ using the elementary algebra methods and the $R$-matrix
approach [4]. Such a decomposition of a given matrix (with non commuting
entries) into the products of matrices of the special type (similar to the
Gauss decomposition) is the particular case of the general
factorization problem [12], which can be considered as a cornerstone of many
constructions of the classical as well as quantum inverse scattering
methods. It is
worthwhile to point out that in different contexts  these
decompositions have been
appeared (sometimes in non explicit form) in many papers on the quantum
deformation [8,12-29,74]. Let us comment, without details, some of them.

In a quite general framework of the quantum double construction
the Gauss decomposition was introduced in the paper [13],
where the universal
triangular objects  ${\cal M}^{\pm}$ were defined. Their matrix
representations on one of the factors $(\rho \otimes 1){\cal M}^{\pm}=M^{\pm}$
are the solutions of the FRT-relations
$$ RM_1^{\pm }M_2^{\pm
}=M_2^{\pm }M_1^{\pm }R,\quad M_1^{+}M_2^{-}=M_2^{-}M_1^{+}.$$
with the $R$-matrix $R=(\rho \otimes \rho){\cal R}$ related to the universal
one ${\cal R}$.
Their product
$M^-M^+$ after the unification of diagonal elements
$(M_{ii}^{-}M_{ii}^{+}=A_{ii})$ gives the $q-$matrix $T$.
Usually the new generators, which are defined by the Gauss decomposition,
have simpler
commutation rules (multiplication) but more complicated expressions
for the coproduct.

The representation of the universal $T$-matrix in terms of the
exponential factors, given in [16] is one of the forms of the Gauss
factorization also. Connection of the Gauss decomposition with Borel
one, which is fundamental for the Borel-Weil construction [30] in the
representation theory also has the quantum group analogs [18, 27-29].

In the dual case of quantum algebras the Gauss decomposition for
the $L$-matrices $L^{\pm}$, which was used in [4] as the very defining
relations,
proves to be useful in the consideration [17] of the equivalence of the
$R$-matrix formulation of the quantum affine algebra
$U_q(\widehat{sl}_2)$ [31]
and the so called "new realization" [32]. The structure
and some properties of
the $R$-matrix may be conveniently described by the Gauss
decomposition [24] also.  Moreover the Gauss factorization was successfully
used in the construction of the quantum group valued
coherent states [18] for the compact quantum algebras.

The simpler structure of the commutation rules among the
Gauss decomposition generators of the
quantum groups simplifies [20, 21] the problem of their
$q$-bosonization [34, 36] (that is their realization by the creation and
annihilation operators of the quantum deformed oscillator [37-39]). For the
dual objects -- quantum algebras (or quantum deformation
of universal enveloping algebras of the classical Lie algebras) this problem
was considered e.g. in [40,41] and for $q$-superalgebras in [42,43]. The
$q$-bosonization of the $L$-operators for the quantum affine
algebra $U_q(\widehat{sl}_2)$, related to the Gauss decomposition, was
analyzed in [22], where connection of the "new realization" of this
algebra at the level 1 with the free field realization of the
Zamolodchikov-Faddeev algebra formulated in the framework of the QISM
was established.

To our surprise many relations among the entries of the quantum
matrix $T$ in terms of the Gauss generators become the elegant
analogs of the standard classical results [33]. The distinctions
appeared in these $q$-analogs consist essentially
in the appearance of the additional $q$-factors and/or
$q-$names only, in the replacement of the
determinants and minors on the quantum determinants and $q$-minors,
respectively and so on. This kind of results were found by many authors and
rediscovered numerously (see, for example, [6,9-11,19-21,27-28]).

In this work we restrict our consideration to those properties of the Gauss
factorization of the quantum group generator matrix $T$ corresponding
to the classical series of the Lie groups and supergroups
which can be obtained by the elementary linear algebra methods
and the $R$-matrix formalism using its fundamental representation.
This allows us to avoid many
difficult and/or cumbersome problems, connected, for example with
the invertibility of some generators and $q-$minors,
with the questions of the original Hopf algebra
structure equivalence with
those generated by new generators appeared in the Gauss decomposition
[13,23], or the questions about connection of the new generators of the
quantum group  with the entries of the $L^{\pm}$-matrices in terms of which
dual quantum algebras are described in the $R$-matrix approach [4].
We leave aside the higher (co)representation Gauss decomposition and
the block matrix cases (see, however Sec.5 on the supergroups).

The paper is set up as follows. The Gauss decomposition of the $GL_q(n)$
is given in Sec.2. In the next Sec.3 we consider the quantum orthogonal
and symplectic groups for their defining relations includes
an additional quadratic constraint . An example of $Sp_q(2)$ is
written in Sec.4 in detail.
Some peculiarities of the quantum supergroups ($Z_2$-graded tensor
product and s-determinant) are mentioned in Sec.5.

{\bf 2. Gauss factorization for the quantum group $GL_q(n)$.} It was mentioned
above that the matrix FRT-relation (\ref{1}) encodes quadratic (commutation)
relations between quantum group generators. The $R$-matrix for the quantum
group $GL_q(n)$ ($A_l$ series) is the lower triangular number
$n^2\times n^2$-matrix [3,4]:
\begin{equation}
\label{2.1}
R=\sum_{i\neq j}^ne_{ii}\otimes e_{jj}+q\sum_{i=1}^ne_{ii}\otimes
e_{ii}+\lambda \sum_{i<j}^ne_{ji}\otimes e_{ij}, \end{equation}
where $\lambda =q-q^{-1}, q\neq 1$, and $e_{kl}$ are the standard matrix
units that is the ($n\times n$)-matrices all entries of which are zeros
except the element $(e_{kl})_{kl}=1$. It is convenient to parametrize the
rows and columns of the $R$-matrix, acting in the space ${\bf C}^n \otimes
{\bf C}^n$, by the pair of indices. With a such
parameterization the element of the $R$-matrix (\ref{2.1}), which stays on
the intersection of the row $(m,n)$ and the column $(p,r)$, has the form
\begin{equation}
\label{2.2}
R_{mn,pr}=(1+(q-1)\delta _{mn})\delta_{mp}\delta _{nr}
+\lambda \theta (m-n)\delta _{np}\delta _{mr},
\end{equation}
where  $$\theta (m-n)=\left\{ \begin{array}{ccc} 1, & {\rm if} &
m>n \\ 0, & {\rm if} & m\leq n \end{array} \right. ,$$ and the
FRT-relation can be represented in the elementwise form as
\begin{equation}
\label{2.3}\sum_{j,k=1}^nR_{mi,jk}t_{jp}t_{kr}=
\sum_{j,k=1}^nt_{ik}t_{mj}R_{jk,pr}.  \end{equation}
Due to (\ref{2.2}), the relations (\ref{2.3}) may be
rewritten as follows
\begin{equation}
\label{2.4}
\begin{array}{c}
t_{i\,j}t_{ik}=qt_{ik}t_{i\,j},\quad \,\;t_{ik}t_{l\,j}=t_{l\,j}t_{ik}, \\
t_{ik}t_{lk}=qt_{lk}t_{ik},\;\quad [t_{i\,j},t_{lk}]=\lambda t_{ik}t_{lj},
\end{array}
\end{equation}
where $1\leq j<k\leq n,\;1\leq i<l\leq n$. It follows from these relations
that the algebra $F_q$ has a rich  subalgebra structure. In
particular, if we omit in the $q$-matrix $T=(t_{mn})$ $k$ arbitrary rows
and $k$ arbitrary columns ($k<n$) then we get the $q$-matrix
related to the $GL_q(n-k)$ subalgebra. In other words the quadratic
commutation relations, defining $GL_q(n)$ are uniquely fixed by the
condition that four generators which stays at the corners of each rectangular
"drawn" in $T$ are subject to the simple well known relations of
the quantum group $GL_q(2)$
\begin{equation}
\label{2.4a}
T=\left(
\begin{array}{cc}
a & b \\
c & d
\end{array}
\right) ;\qquad
\begin{array}{ccc}
ab=qba, & ac=qca, & bc=cb, \\
bd=qdb, & cd=qdc, & [a,d]=\lambda bc.
\end{array}
\end{equation}
The {\bf quantum determinant} is given by the expression [4]
\begin{equation}
\label{2.5}
D_q(T)=\det \nolimits_qT=\sum_\sigma (-q)^{l(\sigma )}t_{1\sigma
(1)}t_{2\sigma (2)}\cdot \ldots \cdot t_{n\sigma (n)},
\end{equation}
where  $l(\sigma)$ is the  length of the substitution $\sigma \in {\bf S}_n$.
The $q$-determinant $\det \nolimits_qT$ is a central element for
the quantum group $GL_q(n)$ [4,6]. For the quantum group $GL_q(2)$ the
quantum determinant is equal to
$${\rm det}_qT=ad-qbc=da-q^{-1}bc.$$
The invertibility of the quantum determinant ${\rm det_q}T$ is the necessary
condition for introducing the Hopf algebra structure into the algebra
$GL_q(n)$. In the following we suppose that this condition is fulfilled.

The {\bf Gauss decomposition} is the quantum matrix factorization of the
following form
\begin{equation}
\label{2.15}
T=T_LT_DT_U=T_LT^{(+)}=T^{(-)}T_U,
\end{equation}
where $T_L=(l_{ik})$ is the strictly lower triangular matrix (with
unities at the main diagonal $l_{kk}=1$); $T_D={\rm diag}(A_{kk})$;
$T_U=(u_{ik})$ is the strictly upper triangular matrix ($u_{kk}=1$);
$T^{(+)}:=T_DT_U$ and $T^{(-)}=T_LT_D$. It is an easy exercise to
get for the $GL_q(2)$ (6)
$l_{21}=c/a, u_{12}=b/qa, A_{11}=a, A_{22}=det_qT/a$.

As in the classical case, with the
mutually commuting entries of the matrix $T$, the Gauss algorithm to
construct the triangular factorization (see for example [33]) can be
applied for the quantum matrix as well. To realize the
Gauss algorithm let us introduce strictly lower triangular
$n\times n$-matrix  $W_L$ and the upper
triangular matrix $T^{(+)}$ by the equation:
\begin{equation}
\label{2.6}
W_LT=T^{(+)}.
\end{equation}
Let $w_k=(w_{k,1},w_{k,2},...,w_{k,k-1})$ and $t_k=
(t_{k,1},t_{k,2},...,t_{k,k-1})$ be the parts of the $k$-th rows of the
matrix $W_L$ and $T$, respectively, and $T_{(k)}$ denotes the matrix which
was received from $T$ by omitting the last $(n-k)$ rows and $(n-k)$
columns. Then in view of the upper triangularity of the matrix  $T^{(+)}$
the condition (\ref{2.6}) gives the system of the linear equations
\begin{equation}
\label{2.7}
w_kT_{(k-1)}=-t_k
\end{equation}
on elements of the row $w_k$ with generators of the algebra
$F_q$ as coefficients. If the $q$-matrix $T_{(k)}$ is invertible (that is
${\rm det_q}T_{(k)} \neq 0$), then the solution of this system has
the form $$w_k=-t_kT_{(k-1)}^{\,\,-1}.$$ We note that the needed for the
existence of a such solution condition corresponds completely to the
classical condition that all diagonal minors of the classical matrix are
different from the zero. The elements of the inverse matrix are given  by the
formula [4,6]:
\begin{equation}
\label{2.10}
(T^{\,-1})_{ij}=(-q)^{i-j}(\det \nolimits_qT)^{-1}M_q(i,j),
\end{equation}
where $M_q(j,i)$ is the $q$-minor, that is the quantum determinant of
the matrix which is obtained from $T$ by omitting of the $i$-th row and
$j$-th column. Thus all the elements of the matrix $W_L$ are uniquely
defined. The condition (\ref{2.6}) allows to find all the nonzero
elements of the matrix $T^{(+)}$. In particular the diagonal elements except
$(T^{(+)})_{11}=t_{11}$, have the form
\begin{equation}
\label{2.11}
\begin{array}{c}
(T^{(+)})_{ii}=\sum\limits_{k=1}^i(W_L)_{ik}t_{ki}
=-\sum\limits_{l,k=1}^{i-1}(-q)^{l-k}t_{il}(\det
\nolimits_qT_{(i-1)})^{-1}M_q(l,k)t_{ki}+t_{ii}= \\
(\det \nolimits_qT_{(i-1)})^{-1}\left[
\sum\limits_{l,k=1}^{i-1}(-q)^{l-k+1}t_{il}M_q(l,k)t_{ki}+(\det
\nolimits_qT_{(i-1)})t_{ii}\right] = \\
(\det \nolimits_qT_{(i-1)})^{-1}\left[
\sum\limits_{k=1}^i(-q)^{i-k}M_q(k,i)t_{ki}\right] =(\det
\nolimits_qT_{(i-1)})^{-1}(\det \nolimits_qT_{(i)}).
\end{array}
\end{equation}
To derive this relation we take into account the relations
$$M_q(n,n)t_{nk}=qt_{nk}M_q(n,n),\quad k<n,$$
and the decomposition ($1\leq k \leq n$)
$$\det\nolimits_qT=\sum_{j=1}^n(-q)^{j-k}t_{kj}M_q(k,j)=
\sum_{j=1}^n(-q)^{k-j}M_q(k,j)t_{jk}.$$
These and some other useful formulas are given in [6]. Let us point out that
the
relation
\begin{equation}
\label{2.12a}
(T^{(+)})_{ii}=
(\det\nolimits_qT_{(i-1)})^{-1}(\det \nolimits_qT_{(i)})
\end{equation}
is the direct $q$-analog of the classical formulas [33]. Then it is easy to
see that the diagonal elements $(T^{(+)})_{ii}$ of the matrix $T^{(+)}$ are
mutually commuting. Using the commutation relations between the
elements $t_{ij}$ of the matrix $T$ and their minors
(which as a rule proved in [6, 9]) one can find the commutation
rules between elements $(T^{(+)})_{ij}$ of the matrix $T^{(+)}$. The
matrix $W_L$ being strictly lower triangular can be inverted. The elements
$(T_L)_{ij}$ of the inverse matrix $(T_L)=(W_L)^{-1}$ are
polynomials on the elements of the matrix $W_L$.  From the relation
(\ref{2.6}) the desired decomposition of the  $q$-matrix $T$ follows
\begin{equation}
\label{2.13}
T=T_LT^{(+)}.
\end{equation}

In the same way one can find the matrix $W_U$, such that $T^{(-)}=TW_U$,
where $T^{(-)}$ is the lower triangular matrix. The
decomposition of the matrix $T$ takes now the form ($W_U^{\,-1}=T_U$)
\begin{equation}
\label{2.13a}
T=T^{(-)}T_U.
\end{equation}
If we multiply from the right the relation (\ref{2.13}) by the matrix
$W_U$, $T^{(-)}=TW_U$, then
$$T^{(-)}=TW_U=T_LT^{(+)}W_U=T_L(T^{(+)}T_U^{\,-1})=T_LT_D,$$
where $T_D=T^{(+)}T_U^{\,-1}$ is the diagonal matrix with the elements
given by (\ref{2.12a}). One has also
$T^{(+)}=T_DT_U$. Hence, we get the Gauss factorization (\ref{2.15}) of
the $q$-matrix $T=T_LT_DT_U$.

The considered above factorization procedure allows to
obtain the expressions for the elements of all matrices participated
in the decomposition in terms of the quantum group $GL_q(n)$ generators
and principal minors. Unfortunately such factorization procedure
is rather cumbersome for it require to find the
commutation rules between the elements of the received matrices. For
the quantum group $GL_q(n)$ most of these formulas can be
extracted from [6]. At the same time for quantum groups corresponding to
other classical series explicit formulas for quantum determinants and minors
as well as the commutation rules are much less known.
This difficulty can be avoided using the contraction procedure
which for the quantum group case was considered, for example, in
\cite{44,45,PK}.

In the fundamental representation the Cartan elements $h_i$ of the Lie
algebra $gl(n)$ are represented by the $n\times n$ matrices
$h_i=\frac 1{2n}e_{ii}$. In the process of the quantum deformation they
are not changed. Moreover they retain to be elements of the Lie type
(primitive)
that is their coproduct has the form $\Delta (h_i)\equiv H_i=
h_i\otimes 1+1\otimes h_i.$ Hence, on the Cartan elements the
coproduct is coincide with the opposite one $\Delta (h_i)=\Delta
^{\prime }(h_i)$, where $\Delta ^{\prime }=P\circ \Delta$ and we have
\begin{equation} \label{2.16}
[R,H_i]=0.
\end{equation}

Let us apply to the FRT-relation (\ref{1}) the similarity transformation
with the matrix
$$K_\gamma =\exp (\sum_{i=1}^n\gamma _ih_i)\otimes \exp
(\sum_{i=1}^n\gamma _ih_i)=\exp (\sum_{i=1}^n\gamma _i\Delta (h_i))=
\exp (\sum_{i=1}^n\gamma _iH_i),$$
where the number coefficients $\gamma _i$
are strictly ordered $\gamma _1>\gamma_2>\ldots >\gamma _n>0$. In view of
the relation (\ref{2.16}), such similarity transformation affects only the
matrices $T_l, (l=1, 2)$:
$$T_l\longrightarrow K_\gamma T_lK_\gamma ^{\,-1}={\rm diag}
(e^{\gamma _i})T{\rm diag}(e^{-\gamma _i}), \qquad t_{ij}\longmapsto
t_{ij}e^{\gamma _i-\gamma _j}.$$
Let us introduce two sets
of the new generators $$t_{\;ij}^{(+)}=\left\{ \begin{array}{cc}
t_{ij}e^{\gamma _i-\gamma _j} & i\leq j \\ t_{ij} & i>j \end{array} \right.
;\quad t_{\;ij}^{(-)}=\left\{\begin{array}{cc} t_{ij} & i<j \\ t_{ij}e^
{\gamma _i-\gamma _j} & i\geq j \end{array} \right. .$$ Let $\gamma_i-
\gamma _j=\gamma _{ij}\varepsilon $, where $\gamma _{ij}>0$, if $i<j$ and
$\gamma _{ij}<0$, if $i>j$. When  $\varepsilon \rightarrow \infty
\;(\varepsilon \rightarrow -\infty )$ in the set \{$t_{\;ij}^{(+)}\}\quad
\{(t_{\;ij}^{(-)})$\} all the matrix elements with $i>j\,\,\,(i<j)$
are disappeared. Thus we constructed the homomorphisms of the algebra
$F_q$ into the algebras $F_q^{(\pm)}$, generated by the upper and lower
triangular $q$-matrices. In this case the commutation rules for the new
generators are defined by the same $R$-matrix:
\begin{equation}
\label{2.19}R\widetilde{T}_1^{(\pm)}\widetilde{T}_2^{(\pm)}=
\widetilde{T}_2^{(\pm )}\widetilde{T}_1^{(\pm)}R.\end{equation}

Similar contraction procedure allows to find the homomorphisms
$\widetilde{T}^{(\pm )}\rightarrow \widetilde{A}$ of the quantum
groups generated by the entries of matrices $\widetilde{T}^{(\pm )}$ into the
group
generators of which are the elements of the diagonal matrix $T_D$. The
commutation rules between elements of this group are also determined by the
FRT-relation
\begin{equation} \label{2.20}R\widetilde{A}_1\widetilde{A}_2=
\widetilde{A}_2\widetilde{A} _1R,\qquad (A_{ij}=\delta _{ij}A_{ij}),
\end{equation}
which due to of the structure of the $R$-matrix for
the $GL_q(n)$ is equivalent to the relation
\begin{equation}
\label{2.21}
\widetilde{A}_1\widetilde{A}_2=
\widetilde{A}_2\widetilde{A}_1.
\end{equation}
This means the commutativity
of the diagonal matrix elements and can be easily proved using also
the FRT-relation (\ref{2.3}) and the $R$-matrix
elements (\ref{2.2}).

Let us return to the relation (\ref{2.19}) and consider once more
contraction with the matrices $K_{1\gamma }=\exp (\sum_{i=1}^n
\gamma _ih_i)\otimes 1$ (or $K_{2\gamma }=1\otimes \exp (\sum_{i=1}^n\gamma
_ih_i)$). The similarity transformation with these operators leaves
without changes in the equality (\ref{2.19}) the matrix
$\widetilde{T}_2^{(\pm )}\;(\widetilde{T}_1^{(\pm )})$, but transforms the
$R$-matrix $R\rightarrow R_\gamma$. For example, the contraction
which maps the matrix $T^{(-)}$ into the diagonal ones ( $\varepsilon
\rightarrow \infty$ under the mentioned above type of ordering in the
set $\left\{ \gamma _i\right\} $), the matrix $R$ goes into
diagonal matrix $R_\gamma \rightarrow R_D$ (where $R_D$ is the diagonal
part of the $R$-matrix). Finally we get the relation
\begin{equation}
\label{2.22}
R_D\widetilde{A}_1\widetilde{T}_{\;2}^{(-)}=
\widetilde{T}_{\;2}^{(-)}\widetilde{A}_1R_D.
\end{equation}
In a similar way one gets the $R$-matrix relation for the matrix
$\widetilde{T}^{(+)}$
\begin{equation} \label{2.23}R_D\widetilde{T}_{\;1}^{(+)}\widetilde{A}_2=
\widetilde{A}_2\widetilde{T}_{\;1}^{(+)}R_D.
\end{equation}

Using the transposition operator ${\cal P}\quad ({\cal P}_{ij,kl}=\delta _{il}
\delta _{jk})$, such that
\begin{equation}
\label{2.24}{\cal P}^2={\cal P},\;{\cal P}T_1{\cal P}=T_2,
\quad {\cal P}T_2{\cal P}=T_1,\quad
{\cal P}R_D{\cal P}=R_D,
\end{equation}
we receive the relations (\ref{2.22}-\ref{2.23}) with the replacement
$1\leftrightarrow 2.$. Finally, after the similarity transformation
of the FRT-relation (\ref{1}) with the operators $K_{\gamma 1},
K_{\gamma 2}$, we get the relation
\begin{equation} \label{2.24a}
R_D\widetilde{T}_{\;1}^{(+)}\widetilde{T}_{\;2}^{(-)}=\widetilde{T}_{\;2}^
{(-)}\widetilde{T}_{\;1}^{(+)}R_D.
\end{equation}
It can be shown that
commutation rules of the new generators defined uniquely from the
commutation relations of the initial generators. This means that the
quantum groups, obtained by this contraction procedure, are isomorphic to
the quantum groups, generated by the respective factors of the Gauss
decomposition (\ref{2.15}) for the quantum group $GL_q(n)$:$\quad\widetilde
{T}^{(+)}\simeq T^{(+)}, \widetilde{A}\simeq T_D, \widetilde{T} ^{(-)}\simeq
T^{(-)}$. Thus we see that commutation rules between generators of the
quantum groups $T_U$ and $T_L$ can be obtained from the relations
(\ref{2.19}-\ref{2.24a}). For example, from the relation (\ref{2.19})
for $T^{(+)}$ and the decomposition $T^{(+)}=T_DT_U$ we have
$$RT_1^{(+)}(T_D)_2T_{U2}=T_2^{(+)}(T_D)_1T_{U1}R $$ or, with the help of
the commutation relations (\ref{2.23}), $$RR_D^{-1}T_{D2}T_1^{(+)}R_DT_{U2}
=R_D^{-1}T_{D1}T_1^{(+)}R_DT_{U1}R.$$ From the relation (\ref{2.20}),
and taking into account that $[R_D,R]=0$, we have finally \begin{equation}
\label{2.25}RT_{U1}R_DT_{U2}=T_{U2}R_DT_{U1}R, \end{equation}
which is the variant of the reflection equations considered for
example in [47-50]. Analogous equation \begin{equation} \label{2.25a}
RT_{L1}R_D^{\;-1}T_{L2}=T_{L2}R_D^{\;-1}T_{L1}R \end{equation} can be
obtained for the matrix $T_L$. Because the diagonal matrices $R_D$ and
$T_{Di}\;(i=1,2)$ commutes with each other, from the relations (\ref{2.22})
and  (\ref{2.23}) follows the equality \begin{equation}
\label{2.26}R_DT_{D1}T_{L2}=T_{L2}T_{D1}R_D \end{equation} \begin{equation}
\label{2.27}R_DT_{U1}T_{D2}=T_{D2}T_{U1}R_D \end{equation}
(and similar relation with substitution $1\leftrightarrow 2$). If we
substitute $T^{(\pm)}$ into the relation (\ref{2.24a}) and take into account
relations (\ref{2.26}) and (\ref{2.27}) we received
\begin{equation} \label{2.27a}T_{U1}T_{L2}=T_{L2}T_{U1}.  \end{equation}
Thus all elements of the $q$-matrix $T_U$ commutes with every element of the
$q$-matrix $T_L$.

On the other hand if the relations (\ref{2.19}-\ref{2.27}) are valid then
every decomposition (\ref{2.15}) fulfills the FRT-relation.  For example,
for the matrix decomposition $T=T^{(-)}T_U$ we have
$$ \begin{array}{c} RT_1T_2-T_2T_1R=
RT_1^{(-)}(T_{U1}T_2^{(-)}-T_2^{(-)}R_D^{\,-1}T_{U1}R_D)T_{U2}- \\
T_2^{(-)}(T_{U2}T_1^{(-)}-T_1^{(-)}R_D^{\,-1}T_{U2}R_D)T_{U1}R+ \\
RT_1^{(-)}T_2^{(-)}R_D^{\,-1}(T_{U1}R_DT_{U2}-R^{\,-1}T_{U2}R_DT_{U1}R),
\end{array}
$$
Three expressions in the brackets () are zeros due to (27) and (24).

Thus the relations (\ref{2.21}),(\ref{2.25})-(\ref{2.27a}) give complete
list of the commutation relations for the elements of the  $q$-matrices
$T_{L},T_{D}$ and $T_U$ from the Gauss decomposition of the
original matrix $T$ for the quantum group $GL_q(n)$. The relations
(\ref{2.19}),
(\ref{2.25})-(\ref{2.25a}), and easily obtained from (\ref{2.26}) and
(\ref{2.27}) relations
\begin{equation} \label{2.26-1}R_DT^{(+)}_{1}T_{L2}=
T_{L2}R_DT^{(+)}_{1};
\end{equation}
\begin{equation}\label{2.27-1}
T^{(-)}_{2}{R_D}^{-1}T_{L1}=T_{L1}T^{(-)}_{2}{R_D}^{-1},
\end{equation}
give the full set of commutation rules between elements of the
matrices $\{T_L, T^{(+)}\}$ and $\{T^{(-)}, T_U\}$. This allows to consider
each of this three sets of elements as the new basis of generators for the
quantum group $GL_q(n)$.
The basis \{$T_{L},T_{D},T_U$\} is particularly convenient one. Diagonal
elements of the matrix $T_{D}$ remind the Cartan generators of the
simple Lie algebra. Quantum determinant (\ref{2.5}) is a central element of
the quantum group $GL_q(n)$ and, as in the case of number matrices can be
expressed as the product of diagonal elements
\begin{equation}
\label{2.28-1}
\det \nolimits_qT=\prod_{i=1}^n(T_D)_{ii}.
\end{equation}
Of course
$\det \nolimits_qT$ commutes with all new generators, this can be easily proved
for the RHS using the relations (\ref{2.26}),(\ref{2.27}). For
example, if we commute the element $u_{mp}$ in succession with every factor
in $\prod_{i=1}^n(T_D)_{ii}$ we get due to the relations
(\ref{2.26}),(\ref{2.27}) the following equalities
 $$ u_{mp}\prod_{i=1}^n(T_D)_{ii}=
\prod_{i=1}^n(T_D)_{ii}\frac{R_{pi,pi}}{R_{mi,mi}}u_{mp}=
\prod_{i=1}^n(T_D)_{ii}u_{mp},$$
because $\prod_{i=1}^nR_{pi,pi}=q$ for every $1\leq p\leq n$. Commutation
relations
of the  $T_{D}$ diagonal elements are defined by the diagonal blocks
of the $R$-matrix (3) (cf (26), (27)). The product of these blocks is
proportional to the unit matrix,
hence the product of the $T_{D}$ diagonal elements is central (which is not the
 case
of the multiparametric $R$-matrix).

{\bf 3. Orthogonal and symplectic quantum groups.}
In this section orthogonal and symplectic quantum groups (that is the
groups of the series $B_l, C_l, D_l $) are considered and we shell omit often
their names in what follows. The commutation relations between generators
of these groups are determined by the FRT-relation with the $R$-matrix [4,5]
\begin{equation}
\label{3.3} \begin{array}{c}
R=q\sum\limits^N_{{i=1} \atop {i\neq i^{\prime}}}
e_{ii}\otimes e_{ii}+e_{
\frac{N+1}2,\frac{N+1}2}\otimes e_{\frac{N+1}2,\frac{N+1}2}
+\sum\limits^N_{{i,j=1} \atop {i\neq j,j^{\prime }}}e_{ii}\otimes
e_{jj} \\ +q^{-1}\sum\limits^N_{{i=1} \atop {i\neq i^{\prime}}}
e_{i^{\prime }i^{\prime }}\otimes e_{ii}+\lambda \sum\limits^N_{{i,j=1}
 \atop {i>j}}e_{ij}\otimes e_{ji}-\lambda \sum\limits_{{i,j=1} \atop
{i>j}}^Nq^{\rho _i-\rho _j}\varepsilon _i\varepsilon _je_{ij}\otimes
e_{i^{\prime }j^{\prime }}.
\end{array}
\end{equation}
in which the second term takes place in the $B$ -- series only. One has in (32)
$N=2n$ for the groups of $C$ - and $D$ - series,
$\quad N=2n+1$
for the $B$ -- series; ${\epsilon}_i=1$ for $B$ -- and $D$ -- series,
$\varepsilon _i=\left\{
\begin{array}{cc} \;1, & i\leq N/2, \\ -1, & i>N/2. \end{array} \right.$ for
$C$ -- series; $i^{\prime }=N+1-i$; $\{{\rho}_i\}$ is a set of numbers fixed
for each group [4, 5].

In the quantum group definition one introduces the additional condition [4]
\begin{equation}
\label{3.2}TCT^tC^{-1}=1=CT^tC^{-1}T,
\end{equation}
which is a quantum analogs of the known conditions for matrices of the
orthogonal and symplectic Lie groups in the fundamental representation.
In (\ref{3.2})  $T$ is the quantum group $q$-matrix;
$T^t$ is the transposed matrix, $C$ is the number matrix [4, 5]
$$C=C_0q^{\rho},\quad {\rho}={\rm
diag}({\rho}_1,{\rho}_2,...,{\rho}_N),\quad
(C_0)_{ij}={\epsilon}_i{\delta}_{ij'}.$$
The relation (\ref{3.2}) has more transparent form in terms of the group
generators
\begin{equation}
\label{3.12} \begin{array}{c}
\sum\limits^N_{k=1}\varepsilon _kq^{{\rho}_{k^{\prime }}}
t_{ik}t_{jk^{\prime }}=\varepsilon _{j^{\prime }}\delta_{ij^{\prime }}
q^{{\rho}_j},
\end{array} \end{equation} \begin{equation} \label{3.13} \begin{array}{c}
\sum\limits^N_{l=1}\varepsilon _lq^{{\rho}_{l^{\prime }}}t_{li}
t_{l^{\prime }j}=\varepsilon \varepsilon _{i^{\prime }}\delta_{i^{\prime}j}
q^{-{\rho}_i}.
\end{array}
\end{equation}
where ${\rho}_i$ and
$\varepsilon _i$ were defined above, $\varepsilon =1$ for $B$ - and
$D$ - series group, $\varepsilon =-1$ for $C$ - series. Among the additional
conditions (\ref{3.12}),(\ref{3.13}) linear independent from the commutation
relations defined by the FRT-equation (\ref{1}) are only those in which
right hand side is unequal to zero (in other words, "diagonal" part of the
relation (\ref{3.2})).

As for the $GL_q(n)$ case, the function algebras $F_q$ of other
series quantum groups also have rich subalgebra structures.
For example, generators which stand in the $q$-matrix $T$ on
crosses of rows with indices $(i_1,i_2,..,i_k)$ and columns with indices
$(j_1,j_2,..,j_k)$ form $GL_q(k)$ subalgebra in $F_q$ if $k\leq n$ and
among row numbers as well as among column numbers there are no pairs
$(i_m,i_l)$
with $i_m=(i_l)^{\prime }=N+1-{i_l}$).

The Gauss factorization of the orthogonal and symplectic groups can be carried
out by the same method as in Sec.2 for the $GL_q(n)$ . The commutation
relations of the generators are determined by the same
formulae (\ref{2.21})-(\ref{2.23}), (\ref{2.25})-(\ref{2.27a}) in the basis
$(T_L,T_D,T_U)$ and by (\ref{2.19}), (\ref{2.25})-(\ref{2.25a}) in the basis
$(T_L,T^{(+)})$ and $(T^{(-)},T_U)$.
Particular properties of orthogonal and symplectic
$R$-matrices mean, that in  deriving the relation (\ref{2.21}) from
(\ref{2.20}), it is necessary also to take into account that in the
considered case
\begin{equation}
\label{3.1}
(T_D)_{ii}(T_D)_{i'i'}=
(T_D)_{jj}(T_D)_{j'j'}\quad 1\leq i,j\leq n;\quad i':=N+1-i.
\end{equation}
This relation is not a restriction for the additional
conditions (\ref{3.2})require stronger relation
\begin{equation}
\label{3.14}
(T_D)_{ii}(T_D)_{i'i'}=1, \quad 1\leq i\leq N.
\end{equation}

Using the commutation relations one can prove that the additional
conditions for the matrices $T^{(\pm)}$ have the same form (\ref{3.2}) as
for $T$, or in terms of the matrix elements (\ref{3.12}),(\ref{3.13}).
Substituting, for example, $T^{(-)}=T_LT_D$ in these
expressions one can obtain the additional conditions for the generators
$(T_L)_{ij}$:
\begin{equation}
\label{3.15}
\sum\limits_{k=1}^N
\varepsilon_kq^{{\rho}_{k^{\prime}}}{\frac{R_{kk^{\prime};kk^{\prime}}}
{R_{kn;kn}}}(T_L)_{mk}(T_L)_{nk^{\prime}}=\varepsilon_{n^{\prime}}
\delta_{mn^{\prime}}q^{{\rho}_n};
\end{equation}
\begin{equation}\label{3.16}
\sum\limits_{l=1}^N\varepsilon_lq^{{\rho}_{l^{\prime}}}{\frac
{R_{mn;mn}}{R_{ml^{\prime};ml^{\prime}}}}(T_L)_{lm}(T_L)_{l^{\prime}n}=
\varepsilon\varepsilon_{m^{\prime}}\delta_{mn^{\prime}}q^{-{\rho}_m}.
\end{equation}
where $(T_L)_{ii}=1;\,\,(T_L)_{ij}=0,$ at $i<j$. The same
relations are valid for $T_U$ - matrix elements. It is essential that the
sums on the left hand sides of (\ref{3.15}),(\ref{3.16}) include linear by
generators terms. It allows us to exclude dependent generators from a
generator list. Finally the number of generators is equal to the dimension of
the
related Lie group. We shell illustrate independent generator choice in
the next section considering the quantum group $Sp_q(2)$ as an example.

{\bf 4. An example: symplectic quantum group} $Sp_q(2)$. In this Sec.
some details of the Gauss factorization of symplectic and orthogonal
quantum groups are considered using $Sp_q(2)$ as an example. One can take
also even simpler case of the quantum group $SO_q(3)$ with the matrix
$(C)_{ij}=(-q)^{i-2}{\delta}_{ij'}, i, j=1, 2, 3 $ in (33). However, due to the
fusion procedure [53] it corresponds to the spin 1 corepresentation of the
$SU_q(2)$ quantum group ( or $GL_q(2)$ if $det_qT$ is not 1 ).
Hence the Gauss decomposition follows from the one
of (6) by tensoring and projecting with the diagonal part
$$( A_{11}, A_{22}, A_{33} )=(a^2, \quad det_qT, \quad (det_qT)^2/a^2), $$
while the additional condition (33) is
a consequence of $T^t{\epsilon}_qT={\epsilon}_q det_qT$ for (6) [61] .  The
entries of $T_L$ and $T_U$ are expressed in terms of two generators
$u_{12}$ and $l_{21}$ ( $[2]_q=q+1/q $)
$$
\begin{array}{c}
l_{31}=l_{21}^2/[2]_q=c^2/a^2;\quad \quad \quad l_{32}=l_{21}/q; \\
u_{13}=u_{12}^2/[2]_q=b^2/q^4a^2;\quad \quad \quad u_{23}=qu_{12}.
\end{array}
$$

The $R$ - matrix for the $Sp_q(2)$ group is
\begin{equation}
\label{5.1}
R=\left(
\begin{array}{cccccccccccccccc}
q & \cdot  & \cdot  & \cdot  & \cdot  & \cdot  & \cdot  & \cdot  & \cdot  &
\cdot  & \cdot  & \cdot  & \cdot  & \cdot  & \cdot  & \cdot  \\
\cdot  & 1 & \cdot  & \cdot  & \cdot  & \cdot  & \cdot  & \cdot  & \cdot  &
\cdot  & \cdot  & \cdot \cdot  & \cdot  & \cdot  & \cdot  & \cdot  \\
\cdot  & \cdot  & 1 & \cdot  & \cdot  & \cdot  & \cdot  & \cdot  & \cdot  &
\cdot  & \cdot  & \cdot  & \cdot  & \cdot  & \cdot  & \cdot  \\
\cdot  & \cdot  & \cdot  & q^{-1} & \cdot  & \cdot  & \cdot  & \cdot  &
\cdot  & \cdot  & \cdot  & \cdot  & \cdot  & \cdot  & \cdot  & \cdot  \\
\cdot  & \lambda  & \cdot  & \cdot  & 1 & \cdot  & \cdot  & \cdot  & \cdot
& \cdot  & \cdot  & \cdot  & \cdot  & \cdot  & \cdot  & \cdot  \\
\cdot  & \cdot  & \cdot  & \cdot  & \cdot  & q & \cdot  & \cdot  & \cdot  &
\cdot  & \cdot  & \cdot  & \cdot  & \cdot  & \cdot  & \cdot  \\
\cdot  & \cdot  & \cdot  & -\lambda /q & \cdot  & \cdot  & q^{-1} & \cdot  &
\cdot  & \cdot  & \cdot  & \cdot  & \cdot  & \cdot  & \cdot  & \cdot  \\
\cdot  & \cdot  & \cdot  & \cdot  & \cdot  & \cdot  & \cdot  & 1 & \cdot  &
\cdot  & \cdot  & \cdot  & \cdot  & \cdot  & \cdot  & \cdot  \\
\cdot  & \cdot  & \lambda  & \cdot  & \cdot  & \cdot  & \cdot  & \cdot  & 1
& \cdot  & \cdot  & \cdot  & \cdot  & \cdot  & \cdot  & \cdot  \\
\cdot  & \cdot  & \cdot  & \lambda /q^3 & \cdot  & \cdot  & \lambda X &
\cdot  & \cdot  & q^{-1} & \cdot  & \cdot  & \cdot  & \cdot  & \cdot  &
\cdot  \\
\cdot  & \cdot  & \cdot  & \cdot  & \cdot  & \cdot  & \cdot  & \cdot  &
\cdot  & \cdot  & q & \cdot  & \cdot  & \cdot  & \cdot  & \cdot  \\
\cdot  & \cdot  & \cdot  & \cdot  & \cdot  & \cdot  & \cdot  & \cdot  &
\cdot  & \cdot  & \cdot  & 1 & \cdot  & \cdot  & \cdot  & \cdot  \\
\cdot  & \cdot  & \cdot  & \lambda Y & \cdot  & \cdot  & \lambda /q^3 &
\cdot  & \cdot  & -\lambda /q & \cdot  & \cdot  & q^{-1} & \cdot  & \cdot  &
\cdot  \\
\cdot  & \cdot  & \cdot  & \cdot  & \cdot  & \cdot  & \cdot  & \lambda  &
\cdot  & \cdot  & \cdot  & \cdot  & \cdot  & 1 & \cdot  & \cdot  \\
\cdot  & \cdot  & \cdot  & \cdot  & \cdot  & \cdot  & \cdot  & \cdot  &
\cdot  & \cdot  & \cdot  & \lambda  & \cdot  & \cdot  & 1 & \cdot  \\
\cdot  & \cdot  & \cdot  & \cdot  & \cdot  & \cdot  & \cdot  & \cdot  &
\cdot  & \cdot  & \cdot  & \cdot  & \cdot  & \cdot  & \cdot  & q
\end{array}
\right)
\end{equation}
where $X=(1+q^{-2});\;Y=(1+q^{-4}).$ This matrix has more nonzero elements
than the same rank quantum group $GL_q(4)$ $R$-matrix,
 in which all the diagonal terms with $q^{-1}$ are changed
to 1, and all the terms with multiplication by
the negative degree of parameter $q$
out of the diagonal are zero. The complicated structure of
the $R$-matrix causes complicated form of the commutation relations.
The last remark refers only to special commutation relations between
elements whose row or/and column indices are connected by the relations
$i=j'$, where $j'=N+1-j$. In the $Sp_q(2)$ case, for which
$N=2n=4$, $j'=5-j$. Other generators have the $GL_q(2)$ commutation
relations as they are elements of the $GL_q(2)$-subalgebras of the quantum
group $Sp_q(2)$. We write down explicitly the "special" relations only:
\begin{equation}
\label{5.2}
\begin{array}{c}
t_{i1}t_{i4}=q^2t_{i4}t_{i1},\quad
t_{i2}t_{i3}=q^2t_{i3}t_{i2}+\lambda t_{i1}t_{i4},\quad 1\leq i\leq4; \\

t_{1i}t_{4i}=q^2t_{4i}t_{1i},\quad
t_{2i}t_{3i}=q^2t_{3i}t_{2i}+\lambda t_{1i}t_{4i},\quad 1\leq i\leq4;  \\

q^{-1}t_{im}t_{i^{\prime }n}=t_{i^{\prime }n}t_{im}+\lambda q^{\rho
_i}\epsilon _i\sum\limits_{j<i}q^{-\rho _j}\epsilon _jt_{jm}t_{j^{\prime
}n},
\quad m\neq n^{\prime};\,\;
\left\{\begin{array}{ccc}
i=1,2; & n<m; & {\rm or} \\
i=3,4; & n>m; &
\end{array} \right. ; \\

q^{-1}t_{im}t_{i^{\prime }n}-\lambda
t_{i^{\prime }m}t_{in}=t_{i^{\prime }n}t_{im}+\lambda q^{\rho _i}\epsilon
_i\sum\limits_{j<i}q^{-\rho _j}\epsilon _jt_{jm}t_{j^{\prime }n}, \\
m\neq n^{\prime},\,\;
i=1,2;\,\,n>m; \\

t_{im}t_{jm^{\prime }}=q^{-1}t_{jm^{\prime }}t_{im}-\lambda
q^{-\rho _m}\epsilon _m\sum\limits_{k>m}q^{\rho _k}\epsilon _kt_{jk^{\prime
}}t_{ik},\quad i\neq j^{\prime};\,\;
\left\{\begin{array}{ccc}
m=1,2; & i>j; & {\rm or} \\
m=3,4; & i<j; &
\end{array} \right. ; \\

t_{im}t_{jm^{\prime }}-\lambda t_{jm}t_{im^{\prime }}=q^{-1}t_{jm^{\prime
}}t_{im}-\lambda q^{-\rho _m}\epsilon _m\sum\limits_{k>m}q^{\rho _k}\epsilon
_kt_{jk^{\prime }}t_{ik},  \\
i\neq j^{\prime },\;m=1,2;\;i<j; \\

q^{-1}t_{ij}t_{i^{\prime }j^{\prime }}-\lambda q^{\rho _i}\epsilon
_i\sum\limits_{k<i}q^{-\rho _k}\epsilon _kt_{kj}t_{k^{\prime }j^{\prime
}}= \qquad\qquad\qquad   \\
\qquad\qquad=q^{-1}t_{i^{\prime }j^{\prime }}t_{ij}-
\lambda q^{-\rho _j}\epsilon
_j\sum\limits_{k>j}q^{\rho _k}\epsilon _kt_{i^{\prime }k^{\prime
}}t_{ik},\,
\left\{\begin{array}{ccc}
j=1,2; & i=3,4; & {\rm or} \\
j=3,4; & i=1,2; &
\end{array} \right. ; \\

q^{-1}t_{ij}t_{i^{\prime }j^{\prime }}-\lambda t_{i^{\prime }j}t_{ij^{\prime
}}-\lambda q^{\rho _i}\epsilon _i\sum\limits_{k<i}q^{-\rho _k}\epsilon
_kt_{kj}t_{k^{\prime }j^{\prime }}= \qquad \qquad \qquad \\
\qquad \qquad \qquad
=q^{-1}t_{i^{\prime }j^{\prime }}t_{ij}-\lambda q^{-\rho _j}\epsilon
_j\sum\limits_{k>j}q^{\rho _k}\epsilon _kt_{i^{\prime }k^{\prime
}}t_{ik},\;\;i,j=1,2;\; \\
\end{array}
\end{equation}
In these formulas
$$
\epsilon_1=\epsilon_2=1,\,\,\epsilon_3=\epsilon_4=-1,\quad
\rho_i=(\rho_1,\rho_2,\rho_3,\rho_4)=(2,1,-1,-2).
$$
Adding the additional relations (\ref{3.12}),(\ref{3.13}) we
get the complete set of commutation relations for the $Sp_q(2)$
quantum group (the FRT-relation together with the additional
constraint (33) yield an overcomplete list).

As in Sec.2 for realization of the classical Gauss algorithm we must
find matrices $W_L$ and $W_U$ solving the
system of equations (\ref{2.7}).  E. g. for $W_L$ one has:
\begin{equation}
\label{5.7}
\begin{array}{c}
w^L_{21}=-t_{21}D_q^{\,-1}[{{1} \atop {1}}];\quad w^L_{31}=(D_q[{{2,3} \atop
{1,3}}]-\lambda D_q[{{1,4} \atop {1,2}}]) D_q^{\,-1}[{{1,2} \atop {1,2}}];
\quad w^L_{32}=-D_q[{{1,3} \atop {1,2}}]D_q^{\,-1}[{{1,2} \atop {1,2}}] ; \\
\hfill \\ w^L_{41}=-q^2t_{41}D_q^{\,-1}[{1 \atop 1}]; \quad
w^L_{42}=-q^2t_{31}D_q^{\,-1}[{1 \atop 1}]; \quad
w^L_{43}=t_{21}D_q^{\,-1}[{1 \atop 1}];
\end{array}
\end{equation}
In these expressions the $q$-minors
$D_q[{{i_1,i_2,...,i_k} \atop {j_1,j_2,...,j_k}}]$ are calculated by the
$GL_q$-rules (\ref{2.5}) for a matrix obtained from $T$ by removing
all the rows and columns except those with indices $(i_1,...,i_k)$
$(j_1,...,j_k)$ respectively. In particular,
$$
\begin{array}{c}
D_q[{{1} \atop {1}}]=t_{11};\qquad D_q[{{1,2} \atop
{1,2}}]=t_{11}t_{22}-qt_{12}t_{21}.
\end{array}
$$
Using (\ref{2.6}) and the additional conditions (\ref{3.12}),(\ref{3.13})
for the matrix $T^{(+)}$ one obtains
$$\begin{array}{c}
T^{(+)}=
\left(
\begin{array}{cccc}
D_q[{{1} \atop {1}}] & t_{12} & t_{13} & t_{14} \\
0 & D_q^{\,-1}[{{1} \atop {1}}]D_q[{{1,2} \atop {1,2}}] &
D_q^{\,-1}[{{1} \atop {1}}]D_q[{{1,2} \atop {1,3}}] &
D_q^{\,-1}[{{1} \atop {1}}](t_{11}t_{24}-q^2t_{14}t_{21}) \\ 0 & 0 &
D_q^{\,-1}[{{1,2} \atop {1,2}}]D_q[{{1} \atop {1}}] &
-D_q^{\,-1}[{{1} \atop {1}}]t_{12} \\ 0 & 0 & 0 &
D_q^{\,-1}[{{1} \atop {1}}]
\end{array}
\right)
\end{array}$$

As in the classical commutative case, the diagonal elements for quantum groups
(see the $GL_q(n)$-case above) are the ratios of diagonal minors
$$
\begin{array}{c}
D_q^{\,-1}[{{1,2,...,k-1} \atop  {1,2,...,k-1}}]
D_q[{{1,2,...,k} \atop {1,2,...,k}}].
\end{array}
$$
In the presented matrix these ratios are simplified
using particular properties of the symplectic quantum groups. Namely,
calculating the $T^{(+)}$ matrix elements on the places
related to the diagonal $GL_q$ - minors the expressions of the
following form appear
\begin{equation}
\label{5.9}
\begin{array}{c}
D^{sp}_q[{{1,2,...,k} \atop
{1,2,...,k}}] =\sum\limits_{\sigma}(-q)^{l(\sigma)}q^{l'(\sigma)}
t_{1,\sigma (1)}t_{2,\sigma (2)}\cdot ...\cdot t_{k,\sigma (k)},
\end{array}
\end{equation}
which have an additional factor $q^{l'(\sigma)}$ in comparison with
(\ref{2.5}).
Here $l'(\sigma)$ is the number of transpositions of "specific" elements
(transposition index). For example, $l'(1,3,2,4)=1$ as 2 and $3=2'$ are
transposed, $l'(1,2,4,3)=0$. The expressions (\ref{5.9}) for $i=3,4$
can be simplified in the following manner
\begin{equation}
\label{5.10}
\begin{array}{c}
D_q^{sp}[{{1,2,3} \atop {1,2,3}}]=D^{sp}_q[{{1} \atop {1}}] \qquad
D^{sp}_q(T)=D^{sp}_q[{{1,2,3,4} \atop {1,2,3,4}}]=1.
\end{array}
\end{equation}
So the formula (\ref{5.9}) can be naturally chosen as a quantum
determinant definition in the symplectic case. This definition
is in agreement with the one given by a geometric consideration [11].

Non zero elements of the matrix $W_L^{\;-1}=T_L$ have the form
(we remind that $T_L=(l_{ij}), T_D=(A_{ii}), T_U=(u_{ij})$)
\begin{equation}
\label{5.11}
\begin{array}{c}
(T_L)_{21}=-w_{21};\quad (T_L)_{32}=-w_{32};\quad (T_L)_{43}=-w_{43}; \\
\hfill \\
(T_L)_{31}=w_{32}w_{21}-w_{31}=t_{31}D_q^{\,-1}[{{1} \atop {1}}]; \\
\hfill \\
(T_L)_{41}=-w_{43}w_{32}w_{21}+w_{43}w_{31}+w_{42}w_{21}-w_{41}=t_{41}
D_q^{\,-1}[{{1} \atop {1}}]; \\
\hfill \\
(T_L)_{42}=w_{43}w_{32}-w_{42}=q^{-1}
D_q[{{1,4} \atop {1,2}}]D_q^{\,-1}[{{1,2} \atop {1,2}}];
\end{array}
\end{equation}
Using the definition $T^{(+)}=T_DT_U$ one can obtain
the matrix elements of $T^{(+)}$.
Similar procedure based on the matrix $W_U$ leads,
naturally, to the same results.

These formulas allow us to calculate commutation relations
among the new basis generators induced by the Gauss decomposition.
Finally, they have the form
\begin{equation}
\label{5.12}
\begin{array}{cc}
{\rm (I)}\quad & [A_{kk},A_{jj}] \hspace{.2cm}
=\hspace{.2cm} 0;\hspace{4.4cm} \\
\hfill \\
{\rm (II)}\quad &
{\left\{\begin{array}{rcl}
l_{21}l_{31} & = & q^2l_{31}l_{21}+q\lambda l_{41}; \\
l_{32}l_{21} & = & q^2l_{21}l_{32}-(q^4-1)l_{31}; \hspace{1cm}\\
l_{31}l_{32} & = & q^2l_{32}l_{31}; \\
\lbrack l_{41},l_{ij}\rbrack & = & 0;
\end{array} \right.}
\\
\hfill \\
{\rm (III)}\quad &
\left\{\begin{array}{rcl}
u_{12}u_{13} & = & q^2u_{13}u_{12}+q\lambda u_{14}; \\
u_{23}u_{12} & = & q^2u_{12}u_{23}-(q^2-q^{-2})u_{13};\hspace{.6cm} \\
u_{13}u_{23} & = & q^2u_{23}u_{13}; \\
\lbrack u_{14},u_{kl}\rbrack & = & 0;
\end{array} \right.
\\
\hfill \\
{\rm (IV)}\quad &
\left\{\begin{array}{rcl}
A_ml_{ij} & = & q^{(\delta _{mj}-\delta _{mj^{\prime }}-\delta
_{mi}+\delta _{mi^{\prime }})}l_{ij}A_m; \\
A_mu_{ij} & = & q^{(\delta
_{im}\,-\delta _{im^{\prime }}-\delta _{jm}+\delta _{jmi^{\prime
}})}u_{ij}A_m;
\end{array} \right.
\\
\hfill \\
{\rm (V)}\quad & [u_{kl},l_{ij}] \hspace{.2cm} =
\hspace{.2cm} 0;\hspace{4cm}
\end{array}
\end{equation}
The additional conditions (33) cause the following connections among the new
generators
\begin{equation}
\label{5.13}
\begin{array}{c}
A_{33}=A_{22}^{\,\,-1};\qquad \quad \qquad A_{44}=A_{11}^{\,\,-1};
\\ l_{42}=q^2l_{31}-l_{21}l_{32};\quad \quad \quad l_{43}=-l_{21}; \\
u_{24}=q^2(u_{13}-u_{12}u_{23});\quad u_{34}=-u_{12};
\end{array}
\end{equation}
Hence, the number of independent generators decreases to 10
i.e. to the dimension of the related classical Lie group.

{\bf 5. Quantum supergroups}. Quantum superalgebras appeared naturally
when quantum inverse scattering method \cite{1, 51, 52} was generalized to
the ${\bf Z}_2$- graded systems \cite{51, 52, 53}.
Related $R$-matrices were considered
in \cite{53, 54, 55, 56, G} and simple examples were presented in
\cite{58}. The works \cite{42, 43, 59}
are devoted to a $q$-bosonization of the $q$-superalgebras.
The simple examples of dual objects, i.e. quantum supergroups, were
investigated in \cite{59}- \cite{PU}.
Necessary definitions
and theorems of the "supermathematics" can be found,
for instance, in \cite{65, 66}.

In this Sec. using the simple examples of the quantum supergroups: $GL_q(1|1)$
and $GL_q(2|1)$ it is shown that with minimal corrections
(referring to the sign factors) most of the
formulas and statements on the Gauss decomposition discussed above are
survived in the supergroup case.

The FRT-relation for the quantum supergroups has the same form as (\ref{1}),
but matrix tensor product includes additional sign factors
(${\pm}$) \cite{51, 52} related to ${\bf Z}_2$-grading \cite{65, 66}. Vector
${\bf Z}_2$ - graded space (superspace) decomposes into the direct sum
of subspaces
$V_0\oplus V_1$ of even and odd vectors on which the parity function
($p(v)=0$ at $v\in V_0$ and $p(w)=1$ at $w\in V_1$) is defined.
As a rule, a vector basis with definite parity $p(v_i)=p(i)=0, 1$ is
using. According to this basis,
the row and column parities are introduced
in the matrix space {\rm End}$(V)$ \cite{65, 66}.
The tensor product of two
even matrices $F, G$ $(p(F_{ij})=p(i)+p(j))$ has the following signs
\cite{51, 52}
\begin{equation} \label{4.1}\left( F\otimes G\right)
_{ij;kl}=(-1)^{p(j)(p(i)+p(k))}F_{ik}G_{jl}. \end{equation}
Due to this prescription $T_2=I\otimes T$ has the same
block-diagonal form as in
the standard (non super) case while $T_1=T\otimes I$ includes the additional
sign factor $(-1)$ for odd elements standing at odd rows of blocks.
The $R$-matrices for supergroups can be extracted, for instance, from
\cite{51}, ...\cite{58}, \cite{G}.
For the $GL_q(n|m)$ quantum supergroup the $R$-matrix
structure is the same as for the $GL_q(n+m)$ but at odd-odd rows $q$ is
changed by $q^{-1}$
\begin{equation} \label{4.2}R=\sum_{i,j}\left( 1-\delta
_{ij}(1-q^{1-2p(i))})\right) e_{ii}\otimes e_{jj}+\lambda
\sum_{i>j}e_{ij}\otimes e_{ji}.
\end{equation}
Let us remind that the tensor
product notation in (49) refers to the graded matrices.
For convenience the
$R$-matrices for $GL_q(1|1)$ and $GL_q(2|1)$ are presented below
$$\begin{array}{cc}
{R^{GL(1|1)}=\left(
\begin{array}{cccc}
q & \cdot  & \cdot  & \cdot  \\
\cdot  & 1 & \cdot  & \cdot  \\
\cdot  & \lambda  & 1 & \cdot  \\
\cdot  & \cdot  & \cdot  & q^{-1}
\end{array}
\right) ;\quad} &
{R^{GL(2|1)}=\left(
\begin{array}{ccccccccc}
q & \cdot  & \cdot  & \cdot  & \cdot  & \cdot  & \cdot  & \cdot  & \cdot  \\
\cdot  & 1 & \cdot  & \cdot  & \cdot  & \cdot  & \cdot  & \cdot  & \cdot  \\
\cdot  & \cdot  & 1 & \cdot  & \cdot  & \cdot  & \cdot  & \cdot  & \cdot  \\
\cdot  & \lambda  & \cdot  & 1 & \cdot  & \cdot  & \cdot  & \cdot  & \cdot
\\
\cdot  & \cdot  & \cdot  & \cdot  & q & \cdot  & \cdot  & \cdot  & \cdot  \\
\cdot  & \cdot  & \cdot  & \cdot  & \cdot  & 1 & \cdot  & \cdot  & \cdot  \\
\cdot  & \cdot  & \lambda  & \cdot  & \cdot  & \cdot  & 1 & \cdot  & \cdot
\\
\cdot  & \cdot  & \cdot  & \cdot  & \cdot  & \lambda  & \cdot  & 1 & \cdot
\\
\cdot  & \cdot  & \cdot  & \cdot  & \cdot  & \cdot  & \cdot  & \cdot  &
q^{-1}
\end{array}
\right) .}
\end{array}$$

The same contraction procedure arguments
as in Sec. 2 result in homomorphisms of
$T=(t_{ij})$ onto triangular $T^{(\pm)}$ and diagonal $T_D$ matrices
as well as in the related $R$-matrix relations. Let us present some of them
pointing out peculiarities of the supergroup case as the main $R$-matrix
properties are the same as those in Secs 2, 3.

{}From the relations
\begin{equation}
\label{4.4}RT_{\quad 1}^{(\pm )}T_{\quad 2}^{(\pm )}=T_{\quad 2}^{(\pm
)}T_{\quad 1}^{(\pm )}R
\end{equation}
due to the $R$-matrix block structure the diagonal element commutation
rules follow
\begin{equation}\label{4.5}
R_D(T_D)_1T_{\quad 2}^{(-)}=T_{\quad 2}^{(-)}(T_D)_1R_D,\quad
R_DT_{\quad 2}^{(+)}(T_D)_1=(T_D)_1T_{\quad 2}^{(+)}R_D.
\end{equation}
For the mutually commutative elements $A_{ii}:\,\,T_D={\rm
diag}(A_{11},A_{22},...)$ one has as above
\begin{equation}
\label{4.6}A_{ii}T^{(\pm )}A_{ii}^{-1}
=(R_D)_{\quad ii}^{\pm 1}T^{(\pm )}(R_D)_{\quad ii}^{\mp 1}.
\end{equation}
However, the diagonal block structure of $R_D$ is
different now. As a consequence the $GL_q(n|m)$ central element is
the ratio of the two products corresponding to the even and odd rows
\begin{equation} \label{4.7}{\rm s-det}_qT
=\left( \frac{\prod_{i=1}^nA_{ii}}{\prod_{k=1}^mA_{n+k n+k}}\right) ,
\end{equation}
which is naturally to call by quantum superdeterminant
($q$-Berezinian). Cases of $R$-matrices depending on spectral parameters
see in \cite{57}.

For the supergroup $GL_q(1|1)$ the commutation relations of the $q$-matrix
elements $T=\left( {a\atop \gamma}{\beta\atop d}\right)$ have the
form \cite{59}-\cite{63}
\begin{equation}\label{4.15}
\begin{array}{ccl}
a\beta=q\beta a, & \beta d=q^{-1}d\beta, & \beta \gamma =-\gamma \beta, \\
a\gamma=q\gamma a, & \gamma d=q^{-1}d\gamma, & ad =da+\lambda\gamma \beta,
\end{array} \,\,\,
{\beta}^2={\gamma}^2=0.
\end{equation}
We use the Greek letters for odd (nilpotent) generators.
The Gauss decomposition generators
\begin{equation}
\label{4.8}T=\left(
\begin{array}{cc}
a & \beta  \\
\gamma  & d
\end{array}
\right) =\left(
\begin{array}{cc}
1 & 0 \\
\varsigma  & 1
\end{array}
\right) \left(
\begin{array}{cc}
A & 0 \\
0 & B
\end{array}
\right) \left(
\begin{array}{cc}
1 & \psi  \\
0 & 1
\end{array}
\right) ,
\end{equation}
are connected with the original ones by the formulas
\begin{equation}
\label{4.9}
A=a,\quad \psi =A^{-1}\beta ,\quad \varsigma =\gamma A^{-1},\quad
B=d-\gamma A^{-1}\beta .
\end{equation}
The relations
\begin{equation}
\label{4.10}
\begin{array}{c}
\lbrack A,B]=0,\quad A\psi =q\psi A,\quad A\varsigma =q\varsigma A,\quad
\psi ^2=\varsigma ^2=0, \\
\psi \varsigma +\varsigma \psi =0,\quad B\psi =q\psi B,\quad B\varsigma
=q\varsigma B
\end{array}
\end{equation}
cause centrality of $GL_q(1|1)$ superdeterminant \cite{59}
\begin{equation} \label{4.11}{\rm s-}\det
\nolimits_qT=AB^{-1}=a^2(ad-q\gamma \beta )^{-1}=a/(d-\gamma a^{-1} \beta ).
\end{equation}

In the $GL_q(2|1)$ case the $q$-matrix of generators has the form
\begin{equation}
\label{4.12}T=
\left(
\begin{array}{ccc}
a & b & \alpha \\
c & d & \beta  \\
\gamma & \delta & f
\end{array}
\right) =
\left(
\begin{array}{cc}
M & \Psi \\
\Phi & f
\end{array}
\right) .
\end{equation}
The even $M$ - matrix elements form $GL_q(2)$ subgroup with the
commutation rules (6) and elements of each ($2\times 2$)
submatrix with even generators at its diagonal form $GL_q(1|1)$
supersubgroup with (\ref{4.15}) - type commutation relations.
The other ones are read as follows
$$ \begin{array}{ccc}
\alpha\beta=-q^{-1}\beta\alpha,\quad & c\alpha=\alpha c,\quad &
b\gamma=\gamma b, \\
\gamma\delta=-q^{-1}\delta\gamma,\quad & d\alpha=\alpha d,\quad &
d\gamma=\gamma d,
\end{array} $$
$$ \begin{array}{cc}
a\beta=\beta a+\lambda c\alpha,\quad &
b\beta=\beta b+\lambda d\alpha, \\
a\delta=\delta a+\lambda \gamma b,\quad &
c\delta=\delta c+\lambda \gamma d.
\end{array} $$

Appearing in the Gauss decomposition
\begin{equation}
\label{4.18}T=
\left(
\begin{array}{ccc}
a & b & \alpha \\
c & d & \beta  \\
\gamma & \delta & f
\end{array}
\right) =
\left(
\begin{array}{ccc}
1 & 0 & 0 \\
u & 1 & 0 \\
v & w & 1
\end{array}
\right) \left(
\begin{array}{ccc}
A & 0 & 0 \\
0 & B & 0 \\
0 & 0 & C
\end{array}
\right) \left(
\begin{array}{ccc}
1 & x & y \\
0 & 1 & z \\
0 & 0 & 1
\end{array}
\right)
\end{equation}
new generators have the following commutation rules
$$\begin{array}{ccl}
Ax=qxA,\quad & Ay=qyA,\quad & Az=zA, \\
Au=quA,\quad & Av=qvA,\quad & Aw=wA, \end{array} $$
$$\begin{array}{ccl}
Bx=q^{-1}xB,\quad & By=yB,\quad & Bz=qzB, \\
Bu=q^{-1}uB,\quad & Bv=vB,\quad & Bw=qwB, \end{array} $$
$$\begin{array}{ccl}
Cx=xC,\quad & Cy=qyC,\quad & Cz=qzC, \\
Cu=uC,\quad & Cv=qvC,\quad & Cw=qwC, \end{array} $$
$$[A,B]=[A,C]=[B,C]=0,\qquad y^2=z^2=v^2=w^2=0$$
$$\begin{array}{ccl}
xy=qyx,\quad & yz=-q^{-1}zy,\quad & qxz-zx=\lambda y, \\
uv=qvu,\quad & vw=-q^{-1}wv,\quad & uw-qw^{-1}=\lambda v, \end{array} $$
$$\begin{array}{c}
[x,u]=[x,v]=[x,w]=0,\quad [u,y]=[u,z]=0, \\
yv+vy=0=yw+wy,\qquad zv+vz=0=zw+wz. \end{array} $$

The superdeterminant
$$ {\rm s-}det_qT=ABC^{-1}= det_qM/C $$
is a central element. The latter expression follows from the block
Gauss decomposition of (59). In particular for the
$GL_q(m|n)$ matrix $T$ in the block form one has (cf. \cite{64})
$$s-det_qT= det_qA/det_q(D-CA^{-1}B)$$
formally the standard expression.

Generalization of the above results to the
$GL_q(m|n)$ and other quantum supergroups looks rather straightforward.
Although,
as usual, the quantum supergroup $OSp_q(1|2)$
 \cite{55, 58} has its own peculiarities. Using notations (60), but with
different grading (0, 1, 0),  for the fundamental corepresentation of $T$ [56]
one has the central
elements $AC=CA=B^2$ and $y=x^2/\omega, v=u^2/\omega, z\simeq x, w\simeq u,
\omega=q^{1/2}-1/q^{1/2}$. The $q$-matrix $T$ has three independent generators
$A, x, u$ while in the undeformed case it is the five parameter supergroup.

{\bf 6. Conclusion}. In this paper we considered the Gauss decomposition
of the quantum groups related to the classical Lie groups and supergroups
by the elementary linear algebra and $R$-matrix methods.
The Gauss factorization yields a new basis for these groups which is
sometimes
more convenient than the original one. Most of
the relations among the Gauss generators are written in the
$R-$matrix form. These commutation
relations are simpler then the original rules. This is especially
clear in the symplectic and orthogonal cases. The additional
conditions completing the $B,C,D$ - series quantum group definition,
allow to extract in terms of the new generators independent ones.
The number of the latter ones is equal to the dimension of the related
classical group. The Gauss factorization leads naturally to
appearance of $q$-analogs of such classical notions as
determinants, superdeterminants and minors. We also
want to stress that the new basis is helpful for study the quantum group
representations. In particular, it was shown in [20, 21], that it
simplifies the $q$-bosonization problem.
As we pointed out in the Introduction it looks that almost any
relation and/or statement on the standard matrices are valid for the
$q-$matrices being appropriately "$q-$deformed". We hope that this
paper contributes to support such a feeling
(cf also \cite{GZ, K, T, P} on the characteristic equations)\footnote
 { After the Russian text of this paper was sent for publication
(Zap. Nauch. Semin. POMI, {\bf 224} (1995) ) we were informed on the Ref.
[74] with detailed analysis of the universal matrix ${\cal T}$ Gauss
decomposition.}.


{\bf Acknowledgement}.The authors would like to thank
L. D. Faddeev, H. Grosse, B. Jur\v co, V. D. Lyakhovsky, S. Z. Pakuliak
and V. O. Tarasov for helpful discussions and for pointing out
important references.

The hospitality and support of the Erwin Schr\"odinger International
Institute for Mathematical Physics during the course of this work is
greatly appreciated.

The research of E.V.D. was partially sponsored by the Russian
Fond of Fundamental Resarches under the grant {\rm N 94-01-01157-a} and
the research of P.P.K. by the DGICYT of Ministerio de Educacion y Ciencia
of Spain.


\begin{center}

\end{center}

\end{document}